\documentclass[12pt]{amsart}

\newcommand{\nexteq}{\displaybreak[0]\\ &=}

\newcommand{\qbinom}[2]{\genfrac{[}{]}{0pt}{}{#1}{#2}}
\newcommand{\GL}{\Gamma\mathrm{L}(V)_H}
\newcommand{\PGL}{P\Gamma\mathrm{L}(V)_H}
\newcommand{\x}{\langle x\rangle}
\newcommand{\cA}{\mathcal{A}}
\newcommand{\cB}{\mathcal{B}}

\newtheorem{thm}{Theorem}

\DeclareMathOperator{\PG}{PG}
\DeclareMathOperator{\GF}{GF}

\title{The twisted Grassmann graph is the block graph of
a design}
\author[A. Munemasa]{Akihiro Munemasa}
\address{Graduate School of Information Sciences,
Tohoku University, Sendai, 980-8579 Japan}
\email{munemasa@math.is.tohoku.ac.jp}

\author[V.D. Tonchev]{Vladimir D. Tonchev}
\address{Department of Mathematical Sciences,
Michigan Technological University, Houghton, MI 49931, USA}
\email{tonchev@mtu.edu} 

\date{July 4, 2009}

\begin{document}
\begin{abstract}
In this note, we show that the twisted Grassmann graph
constructed by van Dam and Koolen is the block graph of
the design constructed by Jungnickel and Tonchev.
We also show that the full automorphism group of the
design is isomorphic to the full automorphism group
of the twisted Grassmann graph.
\end{abstract}
\maketitle

\section{Introduction}
Let $V$ be a $(2e+1)$-dimensional vector space over $\GF(q)$.
If $W$ is a subset of $V$ closed under multiplication by the
elements of $\GF(q)$, then we denote by $[W]$ the set of
$1$-dimensional subspaces (projective points) contained in $W$.
We also denote by 
$\qbinom{W}{k}$ the set of $k$-dimensional subspaces of $W$,
when $W$ is a vector space.
The geometric design $\PG_e(2e,q)$ has $[V]$ as the set of points,
and $\{[W]\mid W\in\qbinom{V}{e+1}\}$ as the set of blocks.
The block graph of this design, where two blocks $[W_1],[W_2]$
are adjacent whenever $\dim W_1\cap W_2=e$, is the Grassmann
graph $J_q(2e+1,e+1)$ which is isomorphic to the Grassmann
graph $J_q(2e+1,e)$.

For each prime power $q$ and an integer $e\ge2$, the twisted
Grassmann graph $\tilde{J}_q(2e+1,e)$
discovered by van Dam and Koolen
is a distance-regular graph with the same
parameters as the Grassmann graph $J_q(2e+1,e)$.
The twisted Grassmann graphs were the first family of
non-vertex-transitive distance-regular graphs with 
unbounded diameter.
We refer the reader to \cite{BI,BCN} for an extensive
discussion of distance-regular graphs, and to
\cite{BFK,FKT} for more information on the twisted
Grassmann graphs.

Jungnickel and the second author \cite{JT} constructed
a family of designs which have the same parameters
as $\PG_e(2e,q)$, and showed that these designs
give the first infinite family of counterexamples
to Hamada's conjecture \cite{H1,H2}. The purpose of this note
is to show that the twisted Grassmann graph is the block
graph of the design constructed in \cite{JT}, just
as the Grassmann graph is the block graph of the
design $\PG_e(2e,q)$.

\section{Statements of the result}
Let $H$ be a fixed hyperplane of $V$.
The twisted Grassmann graph $\tilde{J}_q(2e+1,e)$
(see \cite{vDK}) has a set of vertices
$\cA\cup\cB$, where
\begin{align*}
\cA&=\{W\in\qbinom{V}{e+1}\mid W\not\subset H\},\\
\cB&=\qbinom{H}{e-1}.
\end{align*}
The adjacency is defined as follows:
\[
W_1\sim W_2\iff
\begin{cases}
\dim W_1\cap W_2=e&\text{if }W_1\in\cA,\;W_2\in\cA,\\
W_1\supset W_2&\text{if }W_1\in\cA,\;W_2\in\cB,\\
\dim W_1\cap W_2=e-2&\text{if }W_1\in\cB,\;W_2\in\cB.
\end{cases}
\]

Let $\sigma$ be a polarity of $H$. That is,
$\sigma$ is an inclusion-reversing permutation of
the set of subspaces of $H$, such that $\sigma^2$ is the identity.
Then $\sigma(W_1)\cap\sigma(W_2)=\sigma(W_1+W_2)$ holds
for any subspaces $W_1,W_2$ of $H$. We refer the reader to
\cite{Hir} for details on polarities.

The pseudo-geometric design constructed by Jungnickel and
Tonchev \cite{JT} has $[V]$ as the set of points, and $\cA'\cup\cB'$
as the set of blocks, where
\begin{align*}
\cA'&=\{[\sigma(W\cap H)\cup(W\setminus H)]\mid W\in\cA\},\\
\cB'&=\{[W]\mid W\in\qbinom{H}{e+1}\}.
\end{align*}
It is shown in \cite{JT} that the incidence structure
$([V],\cA'\cup\cB')$ is a $2$-$(v,k,\lambda)$ design, where
\[
v=\frac{q^{2e+1}-1}{q-1},\;
k=\frac{q^{e+1}-1}{q-1},\;
\lambda=\frac{(q^{2e-1}-1)\cdots(q^{e+1}-1)}{(q^{e-1}-1)\cdots(q-1)}.
\]
Moreover, as 
shown in \cite{JT},
the sizes of the intersections of pairs of blocks are
\[
\frac{q^i-1}{q-1}\quad(i=1,\dots,e),
\]
which are exactly the same as those for the geometric design 
$\PG_e(2e,q)$. This leads us to define the block graph
of the design $([V],\cA'\cup\cB')$ in the same manner
as in $\PG_e(2e,q)$, and it turns out that this block
graph is isomorphic to the 
twisted Grassmann graph $\tilde{J}_q(2e+1,e)$.

\begin{thm}\label{thm:main}
The twisted Grassmann graph $\tilde{J}_q(2e+1,e)$
is isomorphic to the block
graph of the design $([V],\cA'\cup\cB')$, where two
blocks are adjacent if and only if their intersection
has size $(q^e-1)/(q-1)$.
\end{thm}
\begin{proof}
We define a mapping $f:\cA\cup\cB\to\cA'\cup\cB'$ by
\[
f(W)=\begin{cases}
[\sigma(W\cap H)\cup(W\setminus H)]&\text{if }W\in\cA,\\
[\sigma(W)]&\text{if }W\in\cB.
\end{cases}
\]
It suffices to show
\begin{equation}\label{eq:1}
W_1\sim W_2\iff
|f(W_1)\cap f(W_2)|=\frac{q^e-1}{q-1}.
\end{equation}
If $W_1,W_2$ are subspaces of $V$, then
\begin{align*}
&\dim\sigma(W_1\cap H)\cap\sigma(W_2\cap H)
\\ &=\dim\sigma(W_1\cap H+W_2\cap H)
%
%
\nexteq
2e-\dim W_1\cap H-\dim W_2\cap H+\dim W_1\cap W_2\cap H
%
%
\nexteq
\begin{cases}
\dim W_1\cap W_2
&\text{if $W_1,W_2\in\cA$, $W_1\cap W_2\subset H$}\\
\dim W_1\cap W_2-1
&\text{if $W_1,W_2\in\cA$, $W_1\cap W_2\not\subset H$}\\
\dim W_1\cap W_2+1
&\text{if $W_1\in\cA$, $W_2\in\cB$,}\\
\dim W_1\cap W_2+2
&\text{if $W_1,W_2\in\cB$}
\end{cases}
\end{align*}

Thus, if $W_1,W_2\in \cA$, then
\begin{align*}
&|f(W_1)\cap f(W_2)|
\\ &=
%
%
|[\sigma(W_1\cap H)\cap\sigma(W_2\cap H)]|+
|[W_1\cap W_2\setminus H]|
\nexteq
\begin{cases}
\displaystyle
\frac{q^{\dim W_1\cap W_2}-1}{q-1}
&\text{if }W_1\cap W_2\subset H,\\[4mm]
\displaystyle
\frac{q^{\dim W_1\cap W_2-1}-1}{q-1}
+\frac{q^{\dim W_1\cap W_2}-q^{\dim W_1\cap W_2-1}}{q-1}
&\text{otherwise}
\end{cases}
\nexteq
\frac{q^{\dim W_1\cap W_2}-1}{q-1},
\end{align*}
and hence (\ref{eq:1}) holds.

Similarly, if $W_1\in\cA$, $W_2\in\cB$, then
\[
|f(W_1)\cap f(W_2)|
=\frac{q^{\dim W_1\cap W_2+1}-1}{q-1},
\]
and hence
\begin{align*}
|f(W_1)\cap f(W_2)|=\frac{q^e-1}{q-1}
&\iff
\dim W_1\cap W_2=\dim W_2
\\ &\iff
W_1\sim W_2.
\end{align*}

Finally, if
$W_1,W_2\in\cB$, then
\[
|f(W_1)\cap f(W_2)|
=
\frac{q^{\dim W_1\cap W_2+2}-1}{q-1}.
\]
and hence (\ref{eq:1}) holds.
\end{proof}

\section{The automorphism group}
Let $\GL$ denote the stabilizer of the hyperplane $H$
in the general semilinear group $\Gamma\mathrm{L}(V)$ on $V$.
For each $\phi\in\GL$, we define a permutation $\phi'$ on
$[V]$ by
\begin{equation}\label{eq:phi'}
\phi'(\x)=\begin{cases}
\sigma\phi\sigma(\x)&\text{if $\x\in[H]$,}\\
\phi(\x)&\text{otherwise,}
\end{cases}
\end{equation}
where $\x$ denotes the $1$-dimensional subspace spanned by
a nonzero element $x\in V$. It is straightforward to verify
that $\phi'$ is an automorphism of the design
$([V],\cA'\cup\cB')$.
Indeed, 
suppose $W\in\qbinom{V}{e+1}$, $W\not\subset H$. Then
\begin{align*}
&\phi'([\sigma(W\cap H)\cup(W\setminus H)])
\\ &=
\{\phi'(\x)\mid \x \in
[\sigma(W\cap H)\cup(W\setminus H)]\}
\nexteq
\{\sigma\phi\sigma(\x)\mid \x \in[\sigma(W\cap H)]
\cup\{\phi(\x)\mid \x \in[W\setminus H]\}
\nexteq
\{\x\mid \sigma\phi(W\cap H)\supset \x\in[H]\}
\cup[\phi(W)\setminus H]
\nexteq
[\sigma(\phi(W)\cap H)\cup\phi(W)\setminus H)]
\\ &\in\cA'.
\end{align*}
Next suppose $W\in\qbinom{H}{e+1}$. Then
\begin{align*}
\phi'([W])&=
\{\sigma\phi\sigma(\x)\mid\x\in[W]\}
\nexteq
\{\x\mid\sigma\phi\sigma(W)\subset\x\in[H]\}
\nexteq
[\sigma\phi\sigma(W)]
\\ &\in\cB'.
\end{align*}
Therefore, $\phi'$ is an automorphism of the design
$([V],\cA'\cup\cB')$.

\begin{thm}\label{thm:aut}
Every automorphism of the design $([V],\cA'\cup\cB')$
is of the form {\rm(\ref{eq:phi'})}, and the full automorphism
group of the design $([V],\cA'\cup\cB')$ is isomorphic to
$\PGL$.
\end{thm}
\begin{proof}
Let $\alpha$ be an automorphism
the design $([V],\cA'\cup\cB')$. By abuse of notation,
denote by the same $\alpha$ the permutation of $\cA'\cup\cB'$
induced by $\alpha$. Then Theorem~\ref{thm:main} implies that
$f^{-1}\alpha f$ is an automorphism of the twisted Grassmann
graph $\tilde{J}_q(2e+1,e)$. Since the automorphism group
of $\tilde{J}_q(2e+1,e)$ is $\PGL$ by \cite{FKT}, there exists
an element $\phi\in\GL$ such that
$f^{-1}\alpha f(W)=f\phi(W)$ for all $W\in\cA\cup\cB$. Then
it is easy to verify that $\alpha(B)=\phi'(B)$ for all
$B\in\cA'\cup\cB'$. 
Indeed,
suppose $W\in\cA$, so that
$[\sigma(W\cap H)\cup(W\setminus H)]\in\cA'$. Then
\begin{align*}
\alpha([\sigma(W\cap H)\cup(W\setminus H)])
&=
\alpha f(W)
\nexteq
f\phi(W)
\nexteq
[\sigma(\phi(W)\cap H)\cup(\phi(W)\setminus H)]
\nexteq
[\sigma\phi\sigma(\sigma(W\cap H))\cup\phi(W\setminus H)]
\nexteq
\phi'([\sigma(W\cap H)\cup(W\setminus H)]).
\end{align*}
Next suppose $W\in\cB$, so that $[\sigma(W)]\in\cB'$.
Then
\begin{align*}
\alpha([\sigma(W)])
&=
\alpha f(W)
\nexteq
f\phi(W)
\nexteq
[\sigma\phi(W)]
\nexteq
[\sigma\phi\sigma(\sigma(W))]
\nexteq
\phi'([\sigma(W)]).
\end{align*}
Therefore $\alpha(B)=\phi'(B)$ for all
$B\in\cA'\cup\cB'$.

Since the action of an automorphism of a
$2$-design on blocks uniquely determines the action on
points if the design has no repeated blocks, 
we obtain the desired result.
\end{proof}



\section*{Acknowledgments}

The second author thanks the Graduate School of Information
Sciences at
Tohoku University, Sendai, Japan,
for the hospitality and support during his visit in June 2009
as a Fulbright Senior Specialist, Project \#3388. This author acknowledges also
partial support by NSA Grant H98230-08-1-0065.

\end{document}